\documentclass[12pt,reqno]{amsart}
\usepackage{amsmath,amssymb,amsfonts,amscd,latexsym,amsthm,mathrsfs}
\usepackage[usenames]{color}
\usepackage[unicode]{hyperref}
\theoremstyle{remark}
\newtheorem{question}{\bf Question}
\begin{document}

\title{Diophantine problems related~to~the~Omega~constant}

\date{23 April 2020}

\author{Wadim Zudilin}
\address{Department of Mathematics, IMAPP, Radboud University, PO Box 9010, 6500~GL Nijmegen, Netherlands}
\email{w.zudilin@math.ru.nl}

\dedicatory{To the memory of Jonathan Sondow (1943--2020)}

\maketitle
%==================================================

The \href{https://en.wikipedia.org/wiki/Omega_constant}{Omega constant}
$$
\Omega=0.567143290409783872999968662210355549753815787186512508135131\dots
$$
is the unique real solution of the transcendental equation
$\Omega e^\Omega=1$, which has several equivalent formulations including $-\log\Omega=\Omega$.
It can be numerically computed via the iterates $x_{n+1}=(1+x_n)/(1+e^{x_n})$ starting, for example, from $x_0=1$.
The constant shows up naturally in relation with Abel's polynomial interpolation \cite[Section~3.5]{Wa20}:
the function $f(x)=e^{\Omega x}$ satisfies $f^{(n)}(n)=1$ for all $n=0,1,2,\dots$, where $f^{(n)}(x)$ denotes the $n$th derivative of $f(x)$.
The transcendence of $\Omega$ directly follows from the Lindemann--Weierstrass theorem (in fact, from its `Lindemann' part!), however the quality of rational (and algebraic) approximations to $\Omega$ remains an open problem.
I learned the following question from J.~Sondow in 2008; it remains completely untouched because of lack of approaches.

\begin{question}
\label{Q1}
Is this true that the Omega constant is \emph{not} a Liouville number?
Give an estimate for the irrationality measure of $\Omega$.
\end{question}

Applying the general lower bounds \cite[Theorem~1]{WN96} for $|e^\theta-\alpha|+|\theta-\beta|$, where $\theta\in\mathbb C\setminus\{0\}$ is fixed and $\alpha,\beta$ are varying algebraic numbers, with the choice $\theta=-\Omega$ and $\alpha=-\beta\in\mathbb Q$ one gets $|\Omega-p/q|>\exp(-C\log^2|q|)$ for some absolute constant $C>0$.
Question~\ref{Q1} explicitly asks whether one can achieve the sharper bound $|\Omega-p/q|>\exp(-C\log|q|)$ for some $C>0$.

\medskip
The Omega constant is closely related to the \href{https://en.wikipedia.org/wiki/Lambert_W_function}{Lambert $W$-function}
which is (a particular branch of) the functional inverse of $ze^z$; in other words, $z=W(ze^z)$.
The `standard' branch of $W(x)$ is real-valued for real $x\ge-1/e$; for it we have $\Omega=W(1)$. 
The branch admits the series expansions
$$
W(x)^k=kx^k\sum_{n=0}^\infty\frac{(n+k)^{n-1}}{n!}\,(-1)^nx^n,
\qquad |x|<\frac1e,
$$
for each $k=1,2,\dots$\,. Again, it follows from Lindemann's theorem that the values of $W(x)$ at non-zero rational (and algebraic) numbers are transcendental,
and a natural question to ask is whether one can give some stronger estimates, than those inplied by \cite{WN96}, for approximations of the values by rational (and algebraic) numbers.

Similar problems may be addressed for the values of $W(x)$ at non-zero algebraic points $\theta$ with $0<|\theta|_p<p^{-1/(p-1)}$ in the $p$-adic settings. Mahler's result \cite{Ma32} from 1932 about the transcendence of $\exp(\theta)$ replaces Lindemann's theorem above. This strategy, however, gives one no way to assign a specific quantity from $\mathbb C_p$ to $W(1)$, perhaps leaving the transcendental number $W(p)\in\mathbb C_p$ as an appropriate candidate for a $p$-adic Omega constant.

Possibly a more intrinsic $p$-adic version of the Lambert $W$-function comes from the formal series inverse $W_p(x)$ of $xE_p(x)$, where
$$
E_p(x)=\exp\biggl(x+\frac{x^p}p+\frac{x^{p^2}}{p^2}+\dots+\frac{x^{p^n}}{p^n}+\dotsb\biggr)
=\prod_{\substack{n=1\\(n,p)=1}}^\infty(1-x^n)^{-\mu(n)/n}\in\mathbb Z_p[[x]]
$$
stands for the \href{https://en.wikipedia.org/wiki/Artin_Hasse_exponential}{Artin--Hasse exponential} (and $\mu(n)$ denotes the M\"obius function). Then $W_p(x)\in\mathbb Z_p[[x]]$ so that the values of the series are well defined inside the $p$-adic unit disk $|x|_p<1$. Observe that
$$
\prod_{n=1}^\infty(1-x^n)^{-\mu(n)/n}=\exp(x)
$$
corresponds to the product side of $E_p(x)$ with $p=1$, and in this sense $W_1(x)$ is exactly the real-valued branch $W(x)$ above.
This setup suggests $\Omega_p=W_p(p)$ as a legitimate $p$-adic Omega constant; alternatively, it can be defined as the limit of the iterates $x_{n+1}=(p+x_n)/(1+E_p(x_n))$ with $x_0=p$.
In absence of Lindemann's theorem for the values of function $E_p(x)$ at algebraic points, already the \emph{irrationality} of $\Omega_p$ looks rather challenging.

\medskip
\noindent
\textbf{Acknowledgements.}
I thank Michel Waldschmidt for related comments and helpful feedback.

%==================================================

\end{document}